\documentclass[a4paper,12pt,reqno]{amsart}
\usepackage{amssymb,amsfonts,amscd,latexsym}
\usepackage{enumerate}
\usepackage{amsmath, mathrsfs,units}
\usepackage{xy}
\xyoption{matrix} \xyoption{arrow}

\numberwithin{equation}{section}

\newcommand{\ad}[1]{\textup{ad}\,{#1}}
\newcommand{\ol}{\overline}

\newcommand{\rest}[2]{{{#1}_{\kern-.5pt|{#2}}}}

\newcommand{\NN}{\mathbb{N}}

\newcommand{\mloc}{M_{\textup{\rm loc}}(A)}
\newcommand{\Mloc}[1]{M_{\textup{\rm loc}}({#1})}
\newcommand{\Mlocbig}[1]{M_{\textup{\rm loc}}\bigl({#1}\bigr)}
\newcommand{\Mlocit}[2]{{M_{\textup{\rm loc}}^{{(#2)}}({#1})}}

\newcommand{\pfi}{\varphi}

\newcommand{\prim}[1]{\textup{\rm Prim}({#1})}
\newcommand{\prima}{{\textup{\rm Prim}(A)}}

\newcommand{\sepa}{{\textup{\rm Sep}(A)}}

\newcommand{\Ice}[1]{{{\mathscr I}_{\mkern-3mu\textup{\rm ce}}({#1})}}
\newcommand{\longrightarrowraised}{{\hbox{\raise.5\jot\hbox{\scriptsize$\longrightarrow$}}}}
\newcommand{\longleftarrowraised}{{\hbox{\raise.5\jot\hbox{\scriptsize$\longleftarrow$}}}}
\newcommand{\dirlim}{{\smash{\underset{\longrightarrowraised}
                   {\operatorname{lim}}}\vphantom{a^{}_{a_f^{}}}}}

\newcommand{\Dirlim}[1]{{\smash{\dirlim{}}\sp{}_{\,{#1}}
                   \vphantom{a^{}_{a_f^{}}}}}
\newcommand\alglim{{\smash{\underset{\longrightarrowraised}
                   {\operatorname{alg\,lim}}}\vphantom{a^{}_{a_f^{}}}}}
\newcommand\Alglim[1]{{\smash{\alglim{}}\sp{}_{\,{#1}}
                   \vphantom{a^{}_{a_f^{}}}}}

\newcommand{\clos}[1]{{\kern.07em{}^c\kern-.11em{#1}}}          
\newcommand{\closa}{{\clos{\kern-.05em{A}}}}                               

\newcommand{\calT}{{\mathcal T}}

\newcommand{\Oprima}{{\mathcal O}_\prima}

\newcommand{\CAlgone}{{\mathcal{C}^*_1}}

\newcommand{\sh}[1]{{\mathfrak{#1}}}
\newcommand{\sha}{{\sh A}}
\newcommand{\shma}{{\sh M_A}}

\newcommand{\shia}{{\sh I_A}}

\newcommand{\shc}{{\sh C}}

\newcommand{\bund}[1]{{\mathsf{#1}}}
\newcommand{\bundA}{{\bund A}}

\newcommand{\sect}[2]{{\Gamma_b({#1},{#2})}}

\def\C*{{\sl C*}-algebra}
\def\Cs*{{\sl C*}-subalgebra}
\def\Cbund*{{\sl C*}-bundle}
\def\Calg*{{\sl C*}-algebraic}
\def\AF/{{\sl A\kern-.5pt F}-algebra}
\def\AW*{{\sl AW*}-algebra}
\def\CXsheaf/{{$\shc(X)$-sheaf}}
\def\CXsheaves/{{$\shc(X)$-sheaves}}

\newtheorem{lem}{Lemma}[section]
\newtheorem{corol}[lem]{Corollary}
\newtheorem{theor}[lem]{Theorem}
\newtheorem{prop}[lem]{Proposition}

\theoremstyle{remark}

\newtheorem{exem}[lem]{\bf Example}

\newtheorem{question*}[]{\bf Question}

\author{Martin Mathieu}
\address{Department of Pure Mathematics\\
Queen's University Belfast\\
Belfast BT7 1NN\\
Northern Ireland}
\email{m.m@qub.ac.uk}

\subjclass[2000]{Primary 46L05. Secondary 46L06, 46M20}
\keywords{Local multiplier algebra, injective envelope, \textsl{C*-}algebra, sheaf theory}

\title[The second local multiplier algebra of a separable \textsl{C*}-algebra]{The second local multiplier algebra of a\\ separable \textsl{C*}-algebra}
\dedicatory{This paper is dedicated to Victor Shulman on his 65th birthday.}

\begin{document}

\begin{abstract}
Several examples of (separable) \C*s with the property that their second (iterated) local multiplier algebra is strictly larger than
the first have been found by various groups of authors over the past few years, thus answering a question originally posed by
G.~K.~Pedersen in 1978. This survey discusses a systematic approach by  P.~Ara and the author to produce such examples on the one hand;
on the other hand, we present new criteria guaranteeing that the second and the first local multiplier algebra of a separable \C* agree.
For this class of \C*s, each derivation of the local multiplier algebra is inner.
\end{abstract}

\maketitle

\section{Introduction}\label{sect:intro}

\noindent
In 1978 the late G.~K.~Pedersen introduced an algebra extension of a general \C* $A$ which he called the ``\C* of essential multipliers''.
His interest in this \C* was stipulated by possible applications in operator theory on \C*s as was the interest of the authors of
the few papers that appeared on this topic around the same time; see~\cite{Ell74} and~\cite{Ell76}.
Having lain dormant for a good while, this theme was taken up again by P.~Ara and the present author who re-discovered the same \C*
independently of Pedersen's work in the late 1980's and coined the terminology ``local multiplier algebra'', introducing, accordingly,
the notation $\mloc$. An important connection with the symmetric algebra of quotients in noncommutative ring theory was
subsequently made, which led to a fuller understanding of the structure of $\mloc$.
For a discussion of this interplay, see~\cite{MM08}.
A comprehensive account can be found in our monograph~\cite{AM03}, which also contains a wealth of applications to a number
of classes of operators between \C*s, thus continuing Pedersen's ideas.

In his seminal paper~\cite{Ped78}, Pedersen asked two questions which spurred quite some research in the past decades.
Iterating the construction of the local multiplier algebra (see the definition in the subsequent section) one obtains the following
tower of \C*s which, a priori, does not have a largest element.
\begin{equation}\label{eq:loc-mult-tower}
A\subseteq\mloc\subseteq\Mloc\mloc\subseteq\ldots
\end{equation}
This led Pedersen to ask
\begin{question*}\label{quest:first}
Is\enspace$\Mloc\mloc=\mloc$\enspace for every \C* $A$?
\end{question*}
The main result in~\cite{Ped78} is the following.
\begin{theor}\label{thm:pedersens}
Let\/ $A$ be a separable \C*. Every derivation\/ $d\colon A\to A$ extends uniquely to a derivation\/
$d\colon\mloc\to\mloc$ and there is\/ $y\in\mloc$ such that\/ $d=\ad y$\/
$(\!$that is, $dx=[x,y]=xy-yx$ for all\/ $x\in\mloc$$)$.
\end{theor}
A derivation of $\mloc$ may or may not leave $A$ invariant; thus Pedersen's second question reads:
\begin{question*}\label{quest:second}
Is every derivation $d\colon\mloc\to\mloc$ inner in~$\mloc$, provided $A$ is a separable \C*?
\end{question*}
It appears that the answer to Question~\ref{quest:second} is still not known in full generality. However,
a negative answer to Question~\ref{quest:first} was provided in~\cite{AM06}. Since then, several other
examples of \C*s with the property that the second local multiplier algebra is strictly larger than the first have
surfaced \cite{AM08}, \cite{ArgFar09}. Some of them are separable \C*s, others are not. On the other hand,
Somerset, in~\cite{Somerset}, gave fairly general conditions on a separable \C*~$A$ implying that both
Questions~\ref{quest:first} and~\ref{quest:second} have positive answers. Until recently, it was not understood
how these two different directions could fit into a more comprehensive framework. This is now achieved in our
paper~\cite{AM11}, and the aim of the present survey is to explain the historical development that led to this
more detailed analysis of the differences between the first and the second local multiplier algebra of a separable \C*.

\section{The History}\label{sect:history}

\noindent
Let $I$ be a closed, two-sided ideal of a \C* $A$ which is essential (that is, $aI=0$ for some $a\in A$ implies that $a=0$).
If $J$ is another such essential ideal of~$A$ which is contained in $I$, then the multiplier algebra $M(I)$ is canonically
embedded as a \Cs* into $M(J)$ by restriction of multipliers to the smaller ideal. In this way, we obtain a directed system
of \C*s with isometric connecting morphisms, where $I$ runs through the directed set $\Ice A$ of all closed, two-sided, essential
ideals of~$A$. The direct limit of this system is $\mloc=\Dirlim{\Ice A}M(I)$, the \textit{local multiplier algebra\/} of~$A$.

There are several other very useful descriptions of the local multiplier algebra, which are all discussed in detail in our monograph~\cite{AM03}.
These lead, e.g., to a representation of the centre $Z=Z(\mloc)$ of $\mloc$ as $Z=\Dirlim{\Ice A}Z(M(I))$
\cite[Section~3.1]{AM03} which cannot be deduced directly from the defining formula for~$\mloc$ above.

Another important characterisation of $\mloc$, which we will rely on heavily in the following, was first obtained by Frank and Paulsen
in~\cite{FrankPaulsen}; see also \cite[Section~4.3]{AM08}. For a \C*~$A$, let us denote by $I(A)$ its \textit{injective envelope\/} as introduced
by Hamana in~\cite{Ham79}; see also~\cite{Paul}. We emphasise that $I(A)$ is not an injective object in the category of \C*s and
*-homomorphisms but in the category of operator spaces and complete contractions. However, it turns out that, nevertheless,
$I(A)$ is a \C* canonically containing $A$ as a \Cs*. For a concise discussion of these facts suited for our purposes, see~\cite{AM08}.

Under this embedding of $A$ into $I(A)$, the local multiplier algebra $\mloc$ is nothing but the completion of the set of all $y\in I(A)$
which act as a multiplier on some $I\in\Ice A$. Since $I(\mloc)=I(A)$ \cite[Proposition~2.14]{AM08}, we see that all higher local multiplier
algebras of $A$ are contained in $I(A)$ wherefore \eqref{eq:loc-mult-tower} improves to
\begin{equation}\label{eq:loc-mult-tower-2}
A\subseteq\mloc\subseteq\Mloc\mloc\subseteq\ldots\subseteq I(A).
\end{equation}
As a result, we can use the injective envelope to study Question~\ref{quest:first}. For instance, if $A$ is commutative, then $\mloc$
is a commutative \AW*, hence injective. It follows that $\mloc=I(A)$ and, since the local multiplier algebra of every \AW* is the algebra
itself \cite[Theorem~2.3.8]{AM03}, we find that $\Mloc\mloc=\Mloc{I(A)}=I(A)=\mloc$, one of the possible ways to affirm Pedersen's question in the
commutative case.

Evidently, $\mloc=M(A)$ for each simple \C* $A$. Since $\Mloc{M(A)}=\mloc$, it follows that $\Mloc\mloc=\mloc$ in this case too.
(We note in passing that $\mloc$ itself can be simple without $A$ being simple unital, see~\cite{AM99}, in which case we also have
$\Mloc\mloc=\mloc$.) On the basis of these two positive answers, it appears to be close at hand to investigate a \C* $A$ which is the tensor product
of a commutative and a simple one:\enspace$A=C\otimes B$ with $C$ commutative (and, without loss of generality, unital) and $B$ simple.
The surprise comes in Section~\ref{sect:results} below where we shall pin down the properties of $C$ and~$B$ determining whether the answer
to Question~\ref{quest:first} is negative or positive.

Another important contribution is due to Somerset who proved in~\cite{Somerset} that the answer to both Question~\ref{quest:first} and
Question~\ref{quest:second} is positive for every unital, separable \C* $A$ which contains sufficiently many maximal ideals; to be precise,
he assumed that the primitive ideal space $\prima$ contains a dense $G_\delta$ subset consisting of closed points. This topological condition
will feature again in Section~\ref{sect:results} below.

The first class of examples of \C*s for which Question~\ref{quest:first} has a negative answer was found in~\cite{AM06}. These are certain unital
separable approximately finite-dimensional \C*s which are primitive and (necessarily) anti-liminal. We employed non-stable \textsl{K}-theory
to describe these \C*s. A very different method, the theory of Hilbert modules over commutative \AW*s, was applied in~\cite{AM08}
to prove that algebras of the form $C(X)\otimes B(H)$, where $H$ is a separable Hilbert space and $X$ is a Stonean space with additional properties,
also provide a negative answer to Question~\ref{quest:first}. From this it is easy to see that some separable, liminal \C*s such as $C[0,1]\otimes K(H)$
have the same property. The latter example was independently found by Argerami, Farenick and Massey in~\cite{ArgFar09}.
Both approaches make use of the injective envelope as well as of formulas for $\mloc$ and $I(A)$ which are (fortunately) available
in these special cases. The same three authors recently studied the local multiplier algebra of certain continuous trace \C*s with similar techniques
\cite{ArgFar10}, \cite{ArgFarNext} leading to the general result that the second local multiplier algebra of those \C*s
is injective \cite[Theorem~6.7]{ArgFarNext}.

\section{The Results}\label{sect:results}

\noindent
In this section, we shall explain how the following dichotomy arises. Let $X$ be a compact metric space which is perfect (that is, contains no
isolated points). Let $B$ be any separable simple \C*. Then we have the following alternative
\begin{equation*}
\xymatrix{&{\kern-.7cm{A=C(X)\otimes B}}\ar@{-}[ddl]_{{B\text{ unital }}}\ar@{-}[ddr]^{{\;B\text{ non-unital}}}& \\
                     & & \\
                     {\Mloc\mloc=\mloc}& &{\kern-1.4cm\Mloc\mloc\neq\mloc}
}
\end{equation*}

\smallskip\noindent
Consequently, there is a plethora of examples of \C*s for which Question~\ref{quest:first} has a negative answer!

\smallskip
The first of the two main results below explains the non-unital case.
\begin{theor}[{\cite[Corollary~3.8]{AM11}}]\label{thm:if-and-only-if}
Let\/ $B$ and\/ $C$ be separable \C*s and suppose that at least one of them is nuclear. Suppose further that\/ $B$ is simple and non-unital
and that\/ $\prim C$ contains a dense $G_\delta$ subset consisting of closed points.
Let\/ $A=C\otimes B$. Then\/ $\mloc=\Mloc\mloc$ if and only if\/ $\prima$ contains a dense subset of isolated points.
\end{theor}
\noindent
The sufficient condition in this theorem applies to any \C*~$A$:
Let $X=\prima$, $X_1$ the set of isolated points in $X$ and $X_2=X\setminus\ol{X_1}$. Then $X_1$ and $X_2$ are disjoint open subsets of~$X$
with corresponding orthogonal closed ideals $I_1=A(X_1)$ and $I_2=A(X_2)$ of~$A$.
If $X_1$ is dense, $I_1$ is the minimal essential closed ideal of $A$ so $\mloc=M(I_1)$. It follows that
\[
\Mloc\mloc=\Mloc{M(I_1)}=\Mloc{I_1}=\mloc.
\]
It is thus surprising that this very general condition is indeed necessary for \C*s of the above type---and therefore provides us with
an easy, systematic way of producing counterexamples. Note that, since $B$ or $C$ is nuclear, the tensor product is unambiguous
and $\prim{C\otimes B}=\prim C\times\prim B$ \cite[Theorem~B.45]{RaeWil}. We shall indicate in Section~\ref{sect:proofs} how
to obtain Theorem~\ref{thm:if-and-only-if}.

The unital case in the above alternative fits into the following more general result which applies to $C(X)\otimes B$ since this
\C* is quasicentral if and only if $B$ is unital (the following formula for the centres applies: $Z(C\otimes B)=Z(C)\otimes Z(B)$
\,\cite[Theorem~3]{Arch75b}).
\begin{theor}[{\cite[Theorem~4.7]{AM11}}]\label{thm:mlocmloc-is-mloc-one}
Let\/ $A$ be a quasicentral separable \C*  such that\/ $\prima$ contains a dense $G_\delta$ subset consisting of closed points.
If\/ $B$ is a \Cs* of\/ $\mloc$ containing\/ $A$ then\/ $\Mloc B\subseteq\mloc$. In particular, $\Mloc\mloc=\mloc$.
\end{theor}
\noindent
Therefore, for every unital separable \C* with Hausdorff primitive spectrum, the first and the second local multiplier algebras
agree; this was already observed in~\cite{Somerset}.

Recall that $A$ is said to be \textit{quasicentral\/} if no primitive ideal of $A$ contains the centre $Z(A)$ of~$A$.
This class of \C*s was introduced and studied initially by Delaroche in the late 1960's \cite{Del67},~\cite{Del68}
and has turned out to be useful on various occasions. Another ``classical'' notion that arises when dealing with topological
spaces associated to a \C* is the one of a separated point. A point $t$ in a topological space $X$ is called \textit{separated\/}
if $t$ and every point $t'$ outside the closure of $\{t\}$ can be separated by disjoint open neighbourhoods.
The set $\sepa$ of all separated points in a separable \C* $A$ is a dense $G_\delta$ subset of $\prima$ and consists precisely of those
$t\in\prima$ for which the norm function $t\mapsto\|a+t\|$ is continuous for every $a\in A$~\cite{Dix68}.
We shall soon make good use of $\sepa$ when we outline the arguments for Theorem~\ref{thm:mlocmloc-is-mloc-one}
in Section~\ref{sect:proofs} below.

But first let us draw an immediate consequence of Theorem~\ref{thm:mlocmloc-is-mloc-one} for Question~2.
\begin{corol}[{\cite[Corollary~4.9]{AM11}}]\label{cor:der-inner-in-mloc}
Let\/ $A$ be a quasicentral separable \C*  such that\/ $\prima$ contains a dense $G_\delta$ subset consisting of closed points.
Then every derivation of\/ $\mloc$ is inner.
\end{corol}
\noindent
The argument is in fact a reduction to Pedersen's theorem, Theorem~\ref{thm:pedersens} above.
Starting with a derivation $d\colon\mloc\to\mloc$ one can construct a separable \Cs* $B$ of $\mloc$ containing $A$ which is $d$-invariant.
By Theorem~\ref{thm:pedersens}, $d$ is inner in $\Mloc B$ which, however, is contained in $\mloc$ by Theorem~\ref{thm:mlocmloc-is-mloc-one}.
For the details of the proof, see~\cite{AM11}. Note that the same argument applies to every separable \C* $A$ with the property
that any separable \Cs* $B$ of $\mloc$ that contains $A$ has its local multiplier algebra $\Mloc B$ contained in~$\mloc$.

Corollary~\ref{cor:der-inner-in-mloc} was obtained in the unital case in~\cite{Somerset}.

\section{Proofs}\label{sect:proofs}

\noindent
We will first attend to the proof of Theorem~\ref{thm:mlocmloc-is-mloc-one}. The key tool are two formulas, valid for an arbitrary \C* $A$,
that describe both $\mloc$ and $I(A)$ in a compatible way. They rest on our sheaf theory for \C*s which was developed in~\cite{AM10}.

Rather than giving the formal definition of  a sheaf of a \C* \cite[Definition~3.1]{AM10}, we shall content ourselves here
with two examples. It is straightforward that both constitute contravariant functors from the category $\Oprima$ of open
subsets of the primitive ideal space of the \C* $A$ with inclusions as arrows into the category $\CAlgone$ of unital \C*s
with unital *-homomorphisms as arrows, that is, \textit{presheaves\/} of unital \C*s over the base space $\Oprima$.
It takes a bit more work to verify the unique gluing property for the multiplier sheaf \cite[Proposition~3.4]{AM10}
whereas it is fairly easy to establish this in the case of the injective envelope sheaf.
\begin{exem}[The multiplier sheaf]\label{exam:multiplier-sheaf}
Let $A$ be a \C* with primitive ideal space $\prima$. We define
\[
\shma\colon\Oprima\to\CAlgone,\quad \shma(U)=M(A(U)),
\]
where $M(A(U))$ denotes the multiplier
algebra of the closed ideal $A(U)$ of $A$ associated to the open subset  $U\subseteq\prima$
and $M(A(U))\to M(A(V))$, $V\subseteq U$, are the restriction homomorphisms.

This is the \textit{multiplier sheaf\/} of $A$ over $\Oprima$.
\end{exem}
\begin{exem}[The injective envelope sheaf]\label{exam:injective-sheaf}
Let $I(B)$ denote the injective envelope of a \C*~$B$. We define
\[
\shia\colon\Oprima\to\CAlgone,\quad\shia(U)=p_UI(A)=I(A(U)),
\]
where $p_U$ denotes the unique central open projection in $I(A)$  such that
$p_U I(A)$ is the injective envelope of~$A(U)$. The mappings
$I(A(U))\to I(A(V))$, $V\subseteq U$, are given by multiplication~by~$p_V$.

This is the \textit{injective envelope sheaf\/} of $A$ over $\Oprima$.
\end{exem}
To every presheaf we can associate in a canonical way an upper semicontinuous \Cbund* \cite[Theorem~5.6]{AM10}.
Let $\sha$ be a presheaf of \C*s over the topological space~$X$. For $t\in X$,  we define the stalk at $t$ by
$\bundA_t:=\dirlim_{t\in U} \sha(U)$ as the direct limit of \C*s of the
directed family $\{ \sha(U)\}$, where $U$ ranges over the family of all open neighbourhoods of $t$ in~$X$.
Take $\bundA:=\bigsqcup_{t\in X} \bundA_t$ and let $\pi(a)=t$ if $a\in\bundA_t$.
For $s\in\sha(U)$ and $t\in U$, we have a canonical *-homomorphism $\sha(U)\to\bundA_t$
and we denote by $s(t)$ the image of $s$ under this mapping. There is a canonical topology on the total
space $\bundA$ (which is uniquely determined in case $X$ is Hausdorff)
turning $(\bundA,\pi ,X)$ into an upper semicontinuous \Cbund* over~$X$.

Let us denote by $\bundA_\shma$ and $\bundA_\shia$, respectively, the \Cbund*s associated in this way to
$\shma$ and $\shia$, respectively.

For any upper semicontinuous \Cbund* $(\bundA,\pi ,X)$ and a subset $Y\subseteq X$, we shall write $\sect Y \bundA$
for the \C* of bounded continuous local sections on~$Y$ \cite[Lemma~5.2]{AM10}.
Moreover, $\calT$ stands for the downwards directed family of dense $G_\delta$ subsets of~$X$.

The following two formulas are the main results, Theorem~7.6 and Theorem~7.7, in~\cite{AM10}.
They require the additional concept of the \textit{derived sheaf\/} of a sheaf of \C*s \cite[Proposition~7.4]{AM10}:
\begin{equation}\label{eq:new-formulas}
\begin{split}
\mloc &=\Alglim{T\in\mathcal T}\,\Gamma_b(T,\bundA_{\shma})\\
            &{}\\
I(A)    &=\Alglim{T\in\mathcal T}\,\Gamma_b(T,\bundA_{\shia}).
\end{split}
\end{equation}
The algebraic direct limit, that is, the limit before completion to a \C*, is already complete since $\prima$ is a Baire space
and hence countable intersections of dense $G_\delta$ subsets remain within~$\calT$.
As $\bundA_\shma$ is a sub-bundle of $\bundA_\shia$, the above formulas~\eqref{eq:new-formulas} open up a new way
to compare the second local multiplier algebra with the first. Indeed, if we take
$y\in\Mloc\mloc\subseteq I(A)$, by~\eqref{eq:new-formulas}, $y$ is contained in some \Cs* $\Gamma_b(T,\bundA_{\shia})$ and
will belong to $\mloc$ once we find $T'\subseteq T$, $T'\in\mathcal T$ such that $y\in\Gamma_b(T',\bundA_{\shma})$.

We shall now outline  how to use this approach to prove Theorem~\ref{thm:mlocmloc-is-mloc-one} for the case $B=\mloc$
(the general case only requires small modifications).
Take $y\in M(J)$ for some closed essential ideal $J$ of $\mloc$. There is $T\in\calT$ such that $y\in\sect{T}{\bundA_\shia}$.
By our assumption on $\prima$ and the remarks before Corollary~\ref{cor:der-inner-in-mloc} we can assume,
without restricting the generality, that $T$ consists of closed separated points of~$\prima$.
By the results in~\cite{AM11}, there is an element $h\in J$ such that $h(t)\ne0$ for all $t\in T$ when viewed as a section on~$T$,
and there is a separable \Cs* $B\subseteq J$  with the properties $AhA\subseteq B$ and $y\in M(B)$.
Now take a countable dense subset $\{b_n\mid n\in\NN\}$  in $B$ and $T_n\in\mathcal T$
such that $b_n\in\sect{T_n}{\bundA_\shma}$ for each $n\in\NN$. In order to simplify the notation, let us put $\bundA=\bundA_\shma$
for the rest of the proof.
Setting $T'=\bigcap_n T_n\cap T\in\mathcal T$, we have $B\subseteq\sect{T'}{\bundA}$ and hence
\[
B_t=\{b(t)\mid b\in B\}\subseteq\bundA_t\qquad(t\in T').
\]
In a next step, we aim to describe the fibres $\bundA_t$ in more detail. In general, there is a homomorphism
$\pfi_t\colon\bundA_t\to\Mloc{A/t}$ \cite[Section~6]{AM10}; however, this need neither be injective nor surjective.
From our hypotheses, we can conclude for each $t\in T'$:
\[
\left.
\begin{matrix}
A &\text{quasicentral ${}\Rightarrow{}$ $A/t$ unital}\\
    & t\enspace\text{closed ${}\Rightarrow{}$ $A/t$ simple}
\end{matrix}
\quad\right\} {}\Rightarrow{}\ \ \Mloc{A/t}=A/t
\]
and, moreover, that $\pfi_t$ is an isomorphism. The surjectivity of $\pfi_t$
rests on the existence of local identities in quasicentral \C*s:
\begin{equation*}
\begin{split}
\forall\ t\in\prima\quad  &\exists\ U_1\subseteq\prima\ \text{open}, t\in U_1,\\
                                                 &\exists\ z\in Z(A)_+,\ \|z\|=1\colon \ z+A(U_2)=1_{A/A(U_2)},
\end{split}
\end{equation*}
where\enspace$U_2=\prima\setminus\ol{U_1}$\enspace\cite[Lemma~4.3]{AM11}.
As a result, $\bundA_t$ is a simple, unital \C* for each $t\in T'$ and thus,
\[
\bundA_t=\bundA_t\,h(t)\,\bundA_t=(A/t)h(t)(A/t)=(AhA)_t\subseteq B_t\subseteq\bundA_t\qquad(t\in T').
\]
Taking $b_t\in B$ with $b_t(t)=1_{\bundA_t}$ we obtain
$y(t)=y(t)\,1_{\bundA_t}=(yb_t)(t)\in\bundA_t$ for all $t\in T'$.
This entails that $y\in\sect{T'}{\bundA_\shma}$  with $T'\subseteq T$, proving that $y\in\mloc$,
as desired.

\smallskip
In the proof of the ``only if''-part in Theorem~\ref{thm:if-and-only-if}, the sheaf theoretic concepts can be pushed to the
background, and we will merely roughly indicate the main argument in the following.
It is easy to reduce the general case to the case when $\prima$ is perfect (i.e., has no isolated points)
using the direct sum decomposition $\mloc=\Mloc{I_1}\oplus\Mloc{I_2}$ (in the notation of the argument
for the ``if''-part directly after the statement of Theorem~\ref{thm:if-and-only-if} above).

The strategy now is the following: to show that $\Mloc\mloc\ne\mloc$ under the hypothesis that $\prima$ is perfect
(and all the other assumptions in Theorem~\ref{thm:if-and-only-if}) we aim to identify a closed essential ideal $K$
of $\mloc$ with the property that $\mloc\varsubsetneq M(K)$. The following \Cs* $K_A$ of $\mloc$ was already introduced
by Somerset in~\cite{Somerset}:\quad$K_A$ is the closure of the set of all elements of the form $\sum_{n\in\NN}a_nz_n$,
where $\{a_n\}\subseteq A$ is a bounded family and $\{z_n\}\subseteq Z(\mloc)$ consists of mutually orthogonal projections.
(These infinite sums exist in $\mloc$ by \cite[Lemma~3.3.6]{AM03}, for example. Note also that $Z(\mloc)$ is countably decomposable since
$A$ is separable.)

For a separable \C* $A$ with the property that $\prima$ contains a dense $G_\delta$ subset of closed points,
the following statements hold (see \cite[Section~2]{AM11}).
\begin{enumerate}[(i)]
\item $K_A$\enspace  is an essential ideal in $\mloc$.
\item If $K_I=K_A$ for all $I\in\Ice A$  then $\Mloc{K_A}=M(K_A)$.
\item Let $y\in I(A)$. If $ya\in K_A$  for all $a\in A$  then $y\in M(K_A)$.
\end{enumerate}
Combining these we obtain the following result.
\begin{prop}[{\cite[Theorem~3.2]{AM11}}]\label{prop:higher-mlocs}\quad
$\Mlocit{A}{3}=\Mlocit{A}{2}=M(K_A)$,
\end{prop}
\noindent
where $\Mlocit{A}{n}=\Mlocbig{\Mlocit{A}{n-1}}$, $n\geq2$ denotes the $n$-fold iterated local multiplier algebra of~$A$.

\smallskip
The main work in the proof of Theorem~\ref{thm:if-and-only-if} consists in constructing an element $q\in M(K_A)\setminus\mloc$.
This uses the special tensor product structure $A=C\otimes B$.
The topological assumptions on $X=\prima=\prim C$ lead us to a dense $G_\delta$ subset $S$ of $X$ consisting of closed separated points
which is a Polish space (using~\cite{Dix68}). Since $S$ itself is perfect,
every non-empty  open subset of $S$ contains an open subset which has non-empty boundary.
This allows us to choose a suitable sequence $(z_n)_{n\in\NN}$ of projections in $Z(M(K_n\otimes B))\subseteq\mloc$, where
$K_n\in\Ice{C}$. As $B$ is simple and non-unital, there is a strictly increasing approximate identity $(e_n)_{n\in\NN}$ of $B$   with
$e_ne_{n+1}=e_n$ and $\|e_{n+1}-e_n\|=1$ for all~$n$.
Put $p_1=e_1$, $p_n=e_n-e_{n-1}$  for $n\geq2$ and set $q_n=\sum_{j=1}^n z_j\otimes p_{2j}$, $n\in\NN$.
We thus obtain an increasing sequence $(q_n)_{n\in\NN}$ in $\mloc_+$  bounded by~$1$.
Since $I(A)$ is monotone complete, $q=\sup_n q_n=\sum_{n=1}^\infty z_n\otimes p_{2n}$ exists
in $I(A)_+$  and has norm~$1$. It remains to show
\begin{itemize}
\item[(a)] $q\in M(K_A)$;
\item[(b)] $q\notin\mloc$.
\end{itemize}
The first assertion is established by using that the bounded central closure $\closa=\ol{AZ}$ is $\sigma$-unital (as $A$ is separable)
together with approximation in the strict topology of $M(K_A)$. For (b), we use special properties of the sequence $(z_n)_{n\in\NN}$ chosen
above relying on the topological properties of~$S$ indicated above. It is shown that, if $q$ were in $\mloc$, this would lead to
a contradiction.

This completes the proof of Theorem~\ref{thm:if-and-only-if}.

\smallskip
It is clear from Proposition~\ref{prop:higher-mlocs}, which holds for every separable \C* $A$ with the property that
$\prima$ contains a dense $G_\delta$ subset of closed points, that the methods discussed above cannot create an example
of a \C* $A$ with the property that $\Mlocit{A}{2}\neq\Mlocit{A}{3}$; in fact, it appears that no such concrete example
is known at present.

\smallskip

\end{document}